\documentclass[12pt,draftcls,onecolumn]{IEEEtran}

\newtheorem{prob}{Problem}
\newtheorem{prop}{Proposition}
\newtheorem{mydef}{Definition}
\newtheorem{lem}{Lemma}

\newtheorem{theorem}{Theorem}

\usepackage{cite}

\ifCLASSINFOpdf
\usepackage[pdftex]{graphicx}
\graphicspath{{../pdf/}{../jpeg/}}
\DeclareGraphicsExtensions{.pdf,.jpeg,.png}
\else
\usepackage[dvipdf]{graphicx}
\graphicspath{{../eps/}}
\DeclareGraphicsExtensions{.eps}
\fi

\usepackage[cmex10]{amsmath}
\usepackage{amssymb}
\usepackage{array}
\usepackage{mdwmath}
\usepackage{mdwtab}
\usepackage{eqparbox}
\usepackage[tight,footnotesize]{subfigure}

\begin{document}

\title{Minimum-Time Quantum Transport with Bounded Trap Velocity}

\author{Dionisis~Stefanatos,~\IEEEmembership{Member,~IEEE,}
        and Jr-Shin~Li,~\IEEEmembership{Member,~IEEE,}
\thanks{D. Stefanatos and J.-S. Li are with the Department
of Electrical and Systems Engineering, Washington University, St. Louis,
MO, 63130 USA, e-mail: dionisis@seas.wustl.edu, jsli@seas.wustl.edu.}
}

\markboth{SUBMITTED TO IEEE TRANSACTIONS ON AUTOMATIC CONTROL}
{Stefanatos \MakeLowercase{\textit{et al.}}: Minimum-Time Quantum Transport with Bounded Trap Velocity}

\maketitle

\begin{abstract}
We formulate the problem of efficient transport of a quantum particle trapped in a harmonic potential which can move with a bounded velocity, as a minimum-time problem on a linear system with bounded input. We completely solve the corresponding optimal control problem and obtain an interesting bang-bang solution. These results are expected to find applications in quantum information processing, where quantum transport between the storage and processing units of a quantum computer is an essential step. They can also be extended to the efficient transport of Bose-Einstein condensates, where the ability to control them is crucial for their potential use as interferometric sensors.

\end{abstract}


\IEEEpeerreviewmaketitle

\section{Introduction}

\IEEEPARstart{D}{uring} the last decades, a wealth of analytical and numerical tools from control theory and optimization have been successfully employed to analyze and control the performance of quantum mechanical systems, advancing quantum technology in areas as diverse as physical chemistry, metrology, and quantum information processing \cite{Mabuchi05}. Although measurement-based feedback control \cite{Wiseman93,Doherty00} and the promising coherent feedback control \cite{Wiseman94,Lloyd00,James08} have gained considerable attention, open-loop control has been proven quite effective since it does not require any quantum measurement, avoiding the associated problems. Controllability results for finite- and infinite-dimensional quantum mechanical systems have been obtained, clarifying the control limits on these systems \cite{Huang83,Li09IEEE,Bloch10}. Some analytical solutions for optimal control problems defined on low-dimensional systems have been derived, yielding novel pulse sequences with unexpected gains compared with those traditionally used \cite{Khaneja01,D'Alessandro01,Boscain02,Sklarz04,Stefanatos04,Stefanatos05,Boscain06, BonnardIEEE,Lapert10,Fisher09,Stefanatos09,Stefanatos10,Stefanatos11,Bonnard10}. Numerical optimization methods, based on gradient algorithms or direct approaches, have also been used intensively to address more complex tasks and to minimize the effect of the ubiquitous experimental imperfections \cite{Peirce88,Khaneja05,Li09,Schulte10,Maximov10,Li11,Ruths11,Doria11}.

At the core of modern quantum technology lies the problem of transfering trapped quantum particles between operational sites by moving the trapping potential. For example, most of the suggested architectures for the implementation of a quantum computer employ the transport of qubits from the storage to the processor unit and back, see \cite{Cirac04}. The transport should be fast and ``faithfull", i.e., the final quantum state should be equivalent to the initial one up to a global phase factor. Ideally, the absence of the vibrational excitations at the final site is required. The high-fidelity transport that satisfies this no-heating condition is characterized as \emph{frictionless}. Note that frictionless quantum transport can be achieved by moving the trapping potential slowly in an adiabatic manner, where the system follows the instantaneous eigenvalues and eigenstates of the time-dependent Hamiltonian. The drawback of this method is the long necessary times which may render it impractical. A way to bypass this problem is to prepare the same final states and energies as with the adiabatic process at a given final time, without necessarily
following the instantaneous eigenstates at each moment. The resulting final state is faithfull while the intermediate states are not.

This nonadiabatic regime, leading to shorter transport times \cite{Masuda10}, provides a privileged area for applying optimal control techniques. Numerical optimization methods have been used to calculate the optimal currents in a segmented Paul trap for fast transport of ions while suppressing vibrational heating \cite{Schulz06,Reichle06}. Fast quantum transport using optical tweezers, where the acceleration is altered in a bang-bang manner, has been demonstrated experimentally \cite{Couvert08} and studied theoretically \cite{TorronteguiPRA11}. For a moving harmonic potential, the limits of faithfull transport with various types of imperfect controls have been evaluated \cite{Murphy09}, while an inverse engineering method using Lewis-Riesenfeld invariants has been employed to achieve efficient quantum transport in short times \cite{TorronteguiPRA11}.

In the present article, we study the problem of minimizing the time of frictionless quantum transport in the case of a harmonic trap moving with a bounded velocity. This is different from the case examined in \cite{Couvert08}, where the acceleration rather than the velocity is bounded. A physical system which can be modeled as a moving potential with bounded speed is the ``magnetic conveyor belt" \cite{Hansel01}. In this system, time-dependent currents in a lithographic conductor pattern create a moving chain of potential wells; atoms are transported in these wells while remaining confined in all three
dimensions. The speed of displacement can be controlled by adjusting the frequencies of the modulating currents. In the next section we formulate the quantum transport with limited trap speed as a time-optimal control problem for a three-dimensional linear system with bounded input. Note that most of the examples presented in the literature are usually limited to two-dimensional systems, which allow the visualization of the optimal synthesis on the plane. The problem is completely solved in section \ref{sec:solution}, where an interesting bang-bang solution is obtained. The present study complements our previous work on minimum-time frictionless cooling of a quantum particle in a harmonic potential \cite{Stefanatos11}.

\section{Optimal Control Formulation of the Quantum Transport Problem}

The evolution of the wavefunction $\psi(x,t)$ of a particle in a
one-dimensional parabolic trapping potential centered around the moving point $s(t)$
is given by the Schr\"{o}dinger equation \cite{Leach77}
\begin{equation}
\label{Schrodinger}
i\hbar\frac{\partial\psi}{\partial t}=\left[-\frac
{\hbar^{2}}{2m}\frac{\partial^{2}}{\partial x^{2}}+\frac{m\omega^{2}}{2}(x-s)^{2}\right] \psi,
\end{equation}
where $m$ is the particle mass and $\hbar$ is Planck's constant; $x$ is a scalar that varies on some compact interval and $\psi$ is a square-integrable function on that interval. We assume that the experimental setup is such that there are essentially no spatial restrictions due to geometrical constraints, for example the system is placed in the middle of a large enough vacuum chamber. When
$s(t)=0$, the above equation can be solved by separation
of variables and the solution is
\begin{equation}
\label{solution_constant}
\psi(x,t)=\sum_{n=0}^{\infty}c_{n} e^{-iE_{n}t/\hbar}\Psi_{n}(x),
\end{equation}
where
\begin{equation}
\label{energy}
E_{n}=\left(n+\frac{1}{2}\right)\hbar\omega,\,n=0,1,\ldots
\end{equation}
are the eigenvalues and
\begin{equation}
\label{eigenstate}
\Psi_{n}(x)=\frac{1}{\sqrt{2^{n}n!}}\left(
\frac{m\omega}{\pi\hbar}\right)  ^{1/4}\exp{\left(  -\frac{m\omega}{2\hbar
}x^{2}\right)  } H_{n}\left(  \sqrt{\frac{m\omega}{\hbar}}x\right)
\end{equation}
are the eigenfunctions of the corresponding time-independent equation
\begin{equation}
\label{time_independent}
\left(-\frac{\hbar^{2}}{2m}\frac{d^{2}
}{d x^{2}}+\frac{m\omega^{2}}{2}x^{2}\right)  \Psi
_{n}=E_{n}\Psi_{n}.
\end{equation}
Here $H_{n}$ in (\ref{eigenstate}) is the Hermite polynomial of degree
$n$. The coefficients $c_{n}$ in (\ref{solution_constant}) can be found from
the initial condition
\begin{equation}
\label{coefficients}
c_{n}=\int_{-\infty}^{\infty}\psi(x,0)\Psi
_{n}(x)dx.\nonumber
\end{equation}

\begin{figure}[t!]
\centering
\includegraphics[width=2.5in]{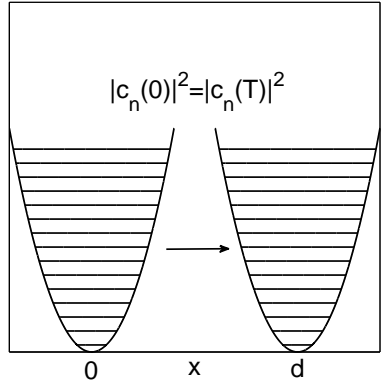}
\caption{Schematic representation of the frictionless atomic transport. The harmonic trapping potential is displaced by $d$, while the populations of all the oscillator levels $n=0,1,2,\ldots$
at the final time are equal to the ones at the initial time. The coefficients $c_n$ should be independent of the spatial coordinate $x$.}
\label{fig:transport}
\end{figure}

Consider now the case where the trap is moving with a bounded velocity $v(t)\in[-V,V], V>0$,
\begin{equation}
\label{mov_trap}
\dot{s}=v(t).
\end{equation}
If $s(0)=0$ and $s(T)=d$, it corresponds to a displacement $d$ of the system in the time interval $[0,T]$, see Fig. \ref{fig:transport}. For
frictionless transport (no vibrational heating), the path $s(t)$ should be
chosen so that the populations of all the oscillator levels $n=0,1,2,\ldots$
for $t=T$ are equal to the ones at $t=0$. In other words, if
\begin{equation}
\label{initial_condition}
\psi(x,0)=\sum_{n=0}^{\infty}c_{n}(0)\Psi_{n}(x),\nonumber
\end{equation}
and
\begin{equation}
\label{final_condition}
\psi(x,T)=\sum_{n=0}^{\infty}c_{n}(T)\Psi_n(x-d),\nonumber
\end{equation}
then frictionless transport is achieved when
\begin{equation}
\label{frictionless_transport}
|c_{n}(T)|^{2}=|c_{n}(0)|^{2},\quad n=0,1,2,\ldots
\end{equation}
This is shown schematically in Fig. \ref{fig:transport}.
We emphasize that the coefficients $c_n$ should be independent of the spatial coordinate $x$.
Among all the paths $s(t)$ that result in (\ref{frictionless_transport}),
we would like to find the one that achieves frictionless transport in minimum time
$T$. In the following we provide a sufficient condition on $s(t)$ for
frictionless transport and we use it to formulate the corresponding time-optimal
control problem.

\begin{prop}
\label{prop:fr_cooling}
If $s(t)$, with $s(0)=0$ and
$s(T)=d$, is such that the
equation
\begin{equation}
\label{moving_oscillator}
\ddot{a}+\omega^{2}(a-s)=0
\end{equation}
has a solution $a(t)$ with $a(0)=0,\dot{a}(0)=0$ and $a(T)=d,\dot{a}(T)=0$, then condition (\ref{frictionless_transport})
for frictionless transport is satisfied.
\end{prop}

\begin{IEEEproof}
Without loss of generality we assume that the initial state is the
eigenfunction corresponding to the $n$-th level $\psi(x,0)=\Psi_n(x)$. We will show that when the hypotheses of Proposition
\ref{prop:fr_cooling} hold then $\psi(x,T)=e^{i\phi_{n}(T)}\Psi_{n}(x-d)$, where $\phi_{n}(T)$ is a global (independent of the spatial
coordinate $x$) phase factor. This and the linearity of (\ref{Schrodinger})
imply that if $\psi(x,0)=\sum_{n=0}^{\infty}c_{n}(0)\Psi_{n}(x)$
then $\psi(x,T)=\sum_{n=0}^{\infty}c_{n}(0)e^{i\phi_{n}(T)}\Psi_{n}(x-d)$, thus condition (\ref{frictionless_transport}) is satisfied.

We follow Leach \cite{Leach77} and consider the ``ansatz"
\begin{equation}
\label{ansatz}
\psi(x,t)=e^{i\left(
\frac{m\dot{a}}{\hbar}x+\phi_n\right)}\Psi_n(x-a),
\end{equation}
where $a(t)$ satisfies (\ref{moving_oscillator})
and the accompanying boundary conditions, while $\phi_n(t)$ is a function of time to be determined, with $\phi_n(0)=0$. Observe that (\ref{ansatz}) corresponds to a wavefunction centered around the moving point $x=a(t)$. The choice of a phase linearly dependent on the spatial coordinate becomes physically transparent if we recall that the momentum operator is $\hat{p}=(\hbar/i)\partial/\partial x$ \cite{Merzbacher98}, so the phase factor in (\ref{ansatz}) gives rise to an average momentum $\langle p\rangle=m\dot{a}$. Note that because of the boundary conditions, we have $\psi(x,0)=\Psi_n(x)$ and $\psi(x,T)=e^{i\phi_{n}(T)}\Psi_{n}(x-d)$ for $t=0$ and $t=T$, respectively, therefore it suffices to show that (\ref{ansatz}) satisfies (\ref{Schrodinger}). Plugging (\ref{ansatz}) into (\ref{Schrodinger}), we obtain
\begin{equation}
\label{intermediate1}
-(m\ddot{a}x+\hbar\dot{\phi}_n)\Psi_n(\chi)=\left[\frac{m\dot{a}^2}{2}-\frac{\hbar^{2}}{2m}\frac{\partial^{2}
}{\partial\chi^{2}}+\frac{m\omega^2}{2}(x-s)^2\right]\Psi_n(\chi),
\end{equation}
where $\chi=x-a$. Since
\begin{equation}
\label{chi}
(x-s)^2=\chi^2+2(a-s)x+s^2-a^2,
\end{equation}
(\ref{intermediate1}) becomes
\begin{equation}
\label{intermediate2}
\left\{m[\ddot{a}+\omega^2(a-s)]x+\hbar\dot{\phi}_n-\frac{\hbar^{2}}{2m}\frac{\partial^{2}
}{\partial\chi^{2}}+\frac{m\omega^2}{2}\chi^2+\frac{m}{2}[\dot{a}^2+\omega^2(s^2-a^2)]\right\}\Psi_n(\chi)=0.
\end{equation}
The coefficient of $x$ is zero because of (\ref{moving_oscillator}). If we additionally use (\ref{time_independent}), then (\ref{intermediate2}) becomes
\begin{equation}
\label{intermediate3}
\left\{\hbar\dot{\phi}_n+E_n+\frac{m}{2}[\dot{a}^2+\omega^2(s^2-a^2)]\right\}\Psi_n(\chi)=0.
\end{equation}
The following choice of $\phi_n$
\begin{equation}
\label{phi}
\phi_n(t)=-\frac{1}{\hbar}\left\{E_nt+\frac{m}{2}\int_0^t[\dot{a}^2+\omega^2(s^2-a^2)]dt\right\}
\end{equation}
assures that (\ref{intermediate3}) is satisfied, so (\ref{ansatz}) is a solution of (\ref{Schrodinger}).
\end{IEEEproof}

We express now the problem of minimum-time frictionless transport using the
language of optimal control. If we set (recall that $V$ is the maximum trap velocity)
\begin{equation}
x_{1}=\frac{\omega}{V}\,a,\quad x_{2}=\frac{\dot{a}}{V},\quad x_{3}=\frac{\omega}{V}\,s,\quad u(t)=\frac{v}{V},
\end{equation}
and rescale time according to $t_{\mbox{new}}=\omega t_{\mbox{old}}$, we
obtain the following linear system with bounded control, equivalent to (\ref{mov_trap}) and (\ref{moving_oscillator}),


\begin{equation}
\label{control_system}
\dot{x}=Ax+u(t)b,
\end{equation}
where now $x=(x_1,x_2,x_3)^T$ and
\begin{equation}
\label{matrices}
A=\left[
\begin{array}{ccc}
0 & 1 & 0\\
-1 & 0 & 1\\
0 & 0 & 0
\end{array}
\right],
\quad\quad
b=
\left[
\begin{array}{c}
0\\0\\1
\end{array}
\right].
\end{equation}
The original transport problem is transformed to the following time-optimal control problem:

\newtheorem{problem}{problem}
\begin{prob}
\label{problem}
Find the control $u(t)$, $|u|\leq 1$, which drives system (\ref{control_system}) from $(0,0,0)$ to $(\gamma,0,\gamma)$, $\gamma=\omega d/V>0$, in minimum time.
\end{prob}

The boundary conditions on $x$ are derived from those on $a$ (see proposition \ref{prop:fr_cooling}) and $s$.



\section{Time-Optimal Solution and Examples}

\label{sec:solution}


Before solving the optimal control problem, we establish the existence and uniqueness of the optimal solution using well known results for linear time-optimal processes. The following two theorems refer to the general linear system
\begin{equation}
\label{gen_system}
\dot{x}=Ax+Bu(t),
\end{equation}
where $x\in\mathbb{R}^n$, $u\in U\subseteq\mathbb{R}^m$, $A\in\mathbb{R}^{n\times n}$, and $B=[b_1|b_2|\ldots|b_m]\in\mathbb{R}^{n\times m}$ ($b_i\in\mathbb{R}^n$).

\begin{theorem}[Controllability of linear systems with bounded controls]
\label{controllability}
Suppose that $A$ is such that all its eigenvalues have real parts equal to zero. Let $U$ be any control set that is a neighborhood of the origin in $\mathbb{R}^m$. Then the linear control system with controls in $U$ is controllable whenever $\bigcup_{k=0}^{n-1}\{A^kb_j, j=1,\ldots,m\}$ spans $\mathbb{R}^n$ (theorem 6, chapter 5 in \cite{Jurdjevic97}).
\end{theorem}

Note that for the single-input case $m=1$, like the system that we study in this article, the above theorem can be directly derived from the null controllability conditions. Recall that the sufficient conditions to be able to bring any initial state of a single-input linear system to zero (null controllability) are that the Kalman matrix has rank $n$, the control $u=0$ belongs to the interior of the control set, and the eigenvalues of matrix $A$ satisfy $\mbox{Re}\,(\lambda_i)\leq 0$ \cite{Lee67}. The full controllability requires additionally the null controllability for the system with matrix $-A$, i.e., $\mbox{Re}\,(\lambda_i)\geq 0$, so that the original system with matrix $A$ can be driven from zero to any final state. The requirements of theorem \ref{controllability}, and especially that for $\mbox{Re}\,(\lambda_i)=0$, are now obvious.

\begin{mydef}[General position condition]
Let the control set $U$ be a convex, closed, and bounded polyhedron in $\mathbb{R}^m$. The matrices $A, B$, and the set $U$ satisfy the general position condition if for every vector $w$, which has the direction of one of the edges of $U$, the vector $Bw$ has the property that it does not belong to any proper subspace of $\mathbb{R}^n$ which is invariant under the operator $A$; i.e., the vectors $Bw, ABw,\ldots, A^{n-1}Bw$ are linearly independent.
\end{mydef}

\begin{theorem}[Existence and uniqueness for linear time-optimal processes]
\label{gen_existence}
Let the control set $U$ be a convex, closed, and bounded polyhedron in $\mathbb{R}^m$ satisfying, along with matrices $A$ and $B$, the general position condition. If there exists at least one control which transfers the state of the system between two points, there also exists a unique optimal control that accomplishes the same transfer (theorems 13 and 11 in \cite{Pontryagin}).
\end{theorem}

\begin{prop}[Existence and uniqueness of the solution for problem \ref{problem}]
\label{existence}
Problem \ref{problem} has a unique optimal solution.
\end{prop}

\begin{IEEEproof}
Matrix $A$ in (\ref{matrices}) has eigenvalues $\pm i, 0$ with zero real parts, while $\mbox{span}\{b, \\Ab, A^2b\}=\mbox{span}\{e_3, e_2, e_1\}=\mathbb{R}^3$, where $e_i$ are the obviously defined unit vectors. Additionally, the control set $U=[-1,1]$ contains the origin. From theorem \ref{controllability} we deduce that system (\ref{control_system}) with $u\in U$ is controllable. So, there exists at least one control which drives the system from the initial to the final point. The general position condition is equivalent to the linear independence of vectors $b, Ab, A^2b$, which is true. From theorem \ref{gen_existence} we conclude that there exists a unique optimal control that accomplishes this transfer.
\end{IEEEproof}

Having established the existence and uniqueness of a solution, we move to solve problem \ref{problem}.
For a constant $\lambda_{0}$ and a row vector $\lambda\in(\mathbb{R}^{3})^\ast$ the control
Hamiltonian for the single-input linear system (\ref{control_system}) is defined as
\[
H=H(\lambda_{0},\lambda,x,u)=\lambda_{0}+\lambda (Ax+ub).
\]
Pontryagin's Maximum Principle \cite{Pontryagin}
provides the following necessary conditions for optimality:

\begin{theorem}[Maximum principle for linear time-optimal processes]
\label{prop:max_principle}
Let $(x_{\ast}(t),u_{\ast}(t))$
be a time-optimal controlled trajectory that transfers the initial condition
$x(0)=x_0$ of system (\ref{control_system}) into the terminal state $x(T)=x_T$. Then it is a necessary
condition for optimality that there exists a constant $\lambda_{0}\leq0$ and
nonzero, absolutely continuous row vector function $\lambda(t)$ such that:

\begin{enumerate}
\item $\lambda$ satisfies the so-called adjoint equation
\[
\dot{\lambda}=-\frac{\partial H}{\partial x}=-\lambda A.
\]

\item For $0\leq t\leq T$ the function $u\mapsto H(\lambda_{0}%
,\lambda(t),x_{\ast}(t),u)$ attains its maximum\ over the control set $U$ at
$u=u_{\ast}(t)$.

\item $H(\lambda_{0},\lambda(t),x_{\ast}(t),u_{\ast}(t))\equiv0$.
\end{enumerate}
\end{theorem}

\begin{mydef}
A control $u:[0,T]\rightarrow[-1,1]$ is said
to be a bang control if $u(t)=+1$ on $[0,T]$ or $u(t)=-1$ on $[0,T]$. A finite concatenation of bang
controls is called a bang-bang control.
\end{mydef}

\begin{prop}
\label{bang_bang}
For problem \ref{problem} extremal controls are bang or bang-bang. The latter controls are $2\pi$-periodic.
\end{prop}

\begin{IEEEproof}
For system (\ref{control_system}) with coefficients given by (\ref{matrices}) we have
\begin{equation}
H(\lambda_0,\lambda,x,u)=\lambda_0+\lambda_1x_2+\lambda_2(x_3-x_1)+\lambda_3u,
\end{equation}
and thus
\begin{align}
\label{l1}\dot{\lambda}_1 &= \lambda_2,\\
\label{l2}\dot{\lambda}_2 &= -\lambda_1,\\
\label{l3}\dot{\lambda}_3 &= -\lambda_2.
\end{align}

Observe that $H$ is a linear function of the bounded control variable $u$. The coefficient of $u$ in $H$ is $\Phi=\lambda_3$, the so-called \emph{switching
function}. According to the maximum principle, point 2 above, the optimal
control is given by $u=\mbox{sign}\,\Phi$, if $\Phi\neq0$. From (\ref{l1}) and (\ref{l3}) we obtain $\lambda_1+\lambda_3=c$, a constant, so $\Phi=\lambda_3=c-\lambda_1$. Also, from (\ref{l1}) and (\ref{l2}) we get $\ddot{\lambda}_1+\lambda_1=0$ (harmonic oscillator), and then $\lambda_1(t)=A\sin(t+\theta)$, where $A$ and $\theta$ are constants. Thus
\begin{equation}
\label{switching_function}
\Phi(t)=c-A\sin(t+\theta).
\end{equation}
The constants $A$ and $c$ cannot be simultaneously equal to zero since $A=0$ implies $\lambda_1=\lambda_2=0$ and $c=\lambda_1=0$ implies $\lambda_3=0$, in contradiction with maximum principle which requires $\lambda=(\lambda_1,\lambda_2,\lambda_3)\neq 0$. Thus the extremal controls are obviously bang or bang-bang, with the latter being $2\pi$-periodic.
\end{IEEEproof}

\begin{figure*}[t]
 \centering
		\begin{tabular}{cc}
     	\subfigure[$\ $First part]{
	            \label{fig:flash1}
	            \includegraphics[width=.45\linewidth]{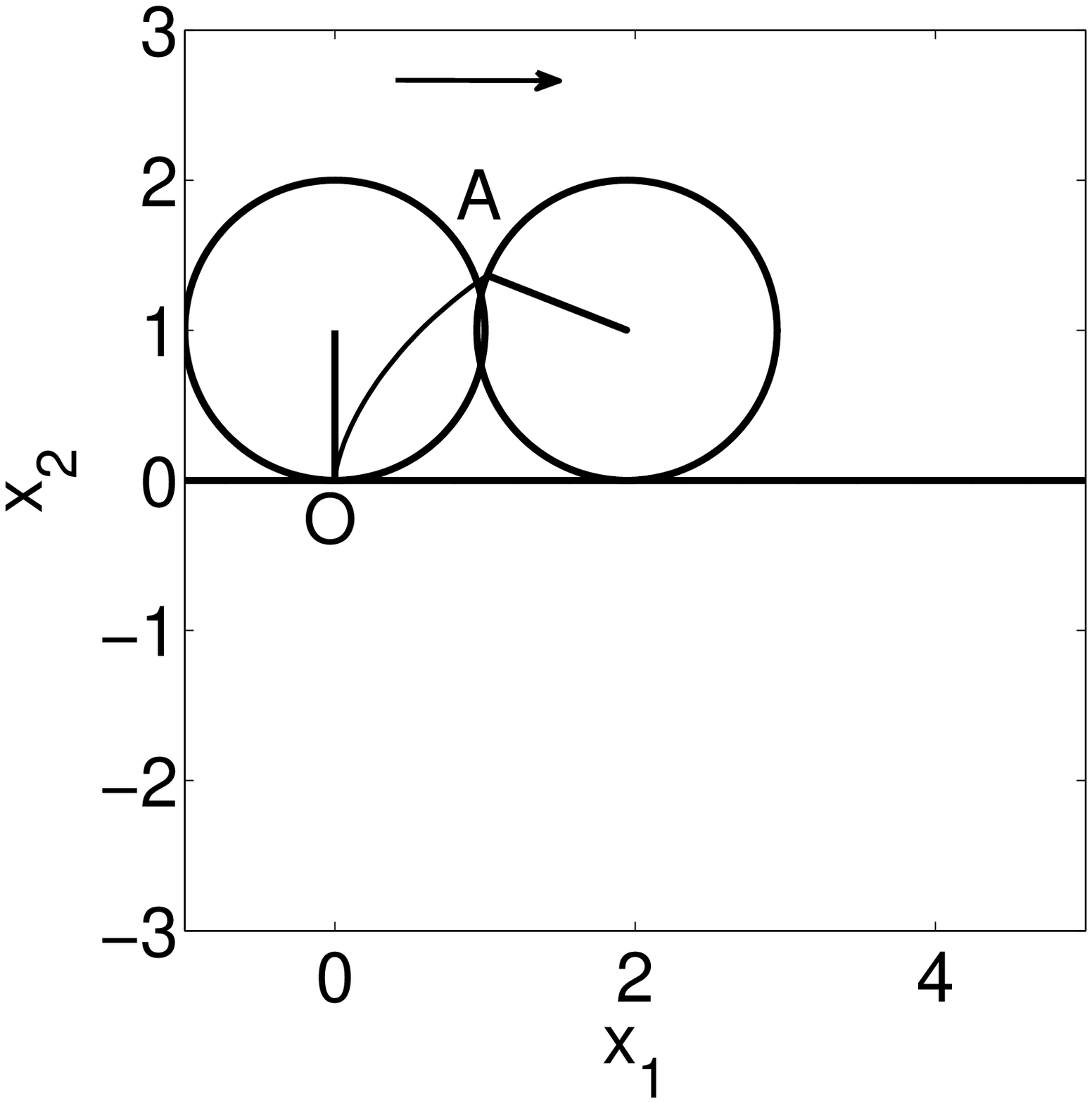}} &
	        \subfigure[$\ $Second part]{
	            \label{fig:flash2}
	            \includegraphics[width=.45\linewidth]{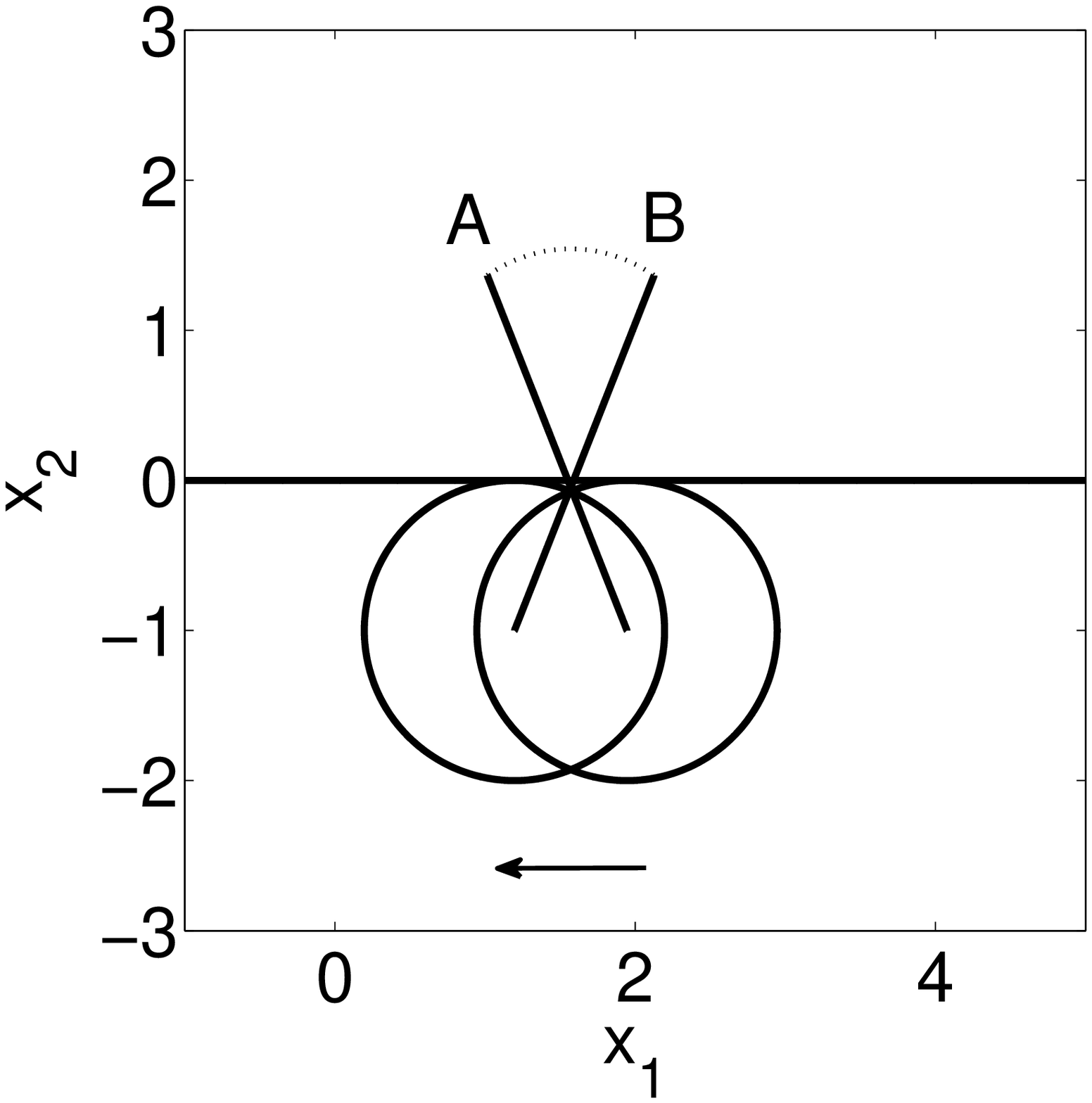}} \\
	        \subfigure[$\ $Third part]{
	            \label{fig:flash3}
	            \includegraphics[width=.45\linewidth]{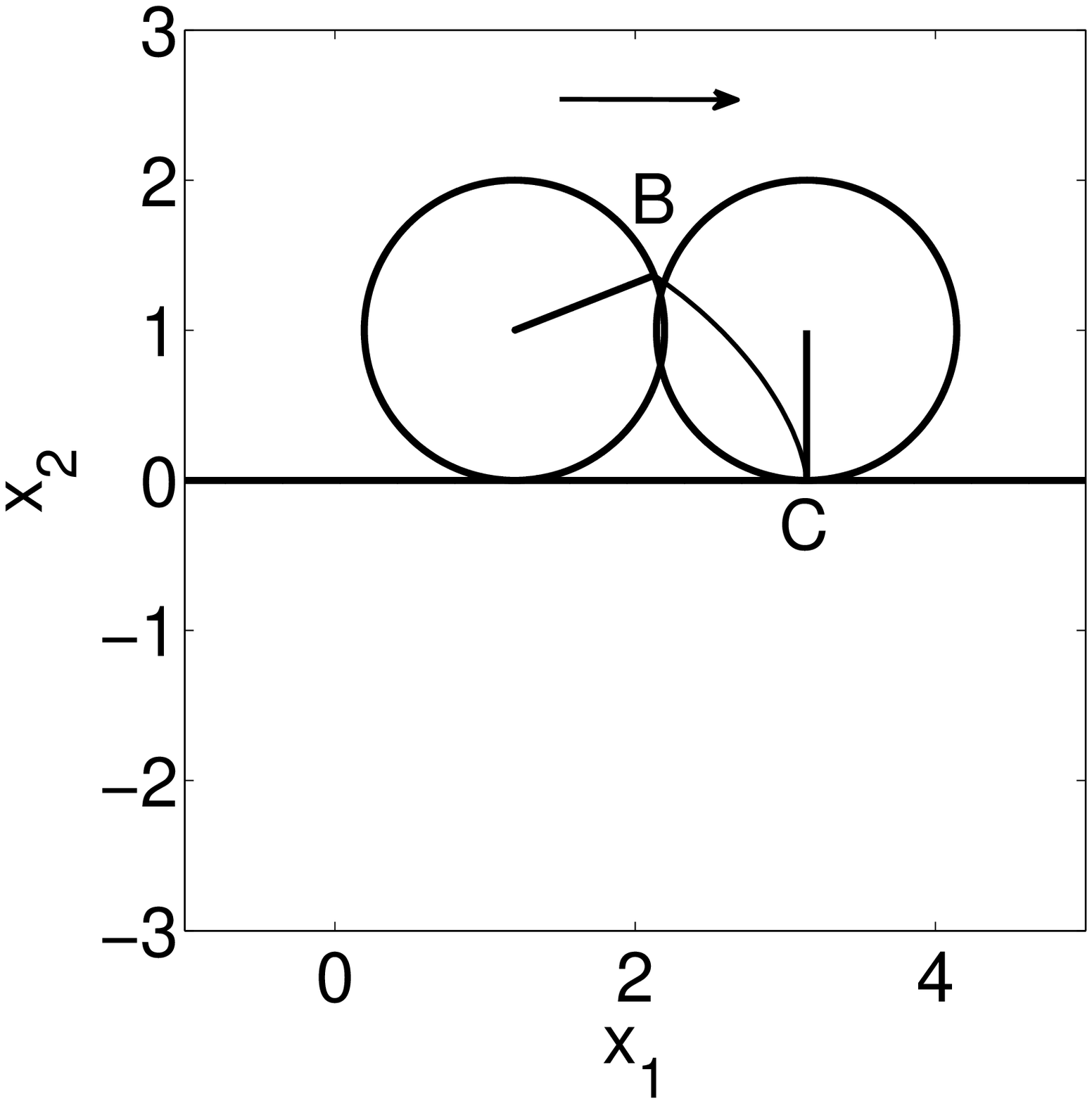}} &
			\subfigure[$\ $Total trajectory]{
	            \label{fig:total_traj}
	            \includegraphics[width=.45\linewidth]{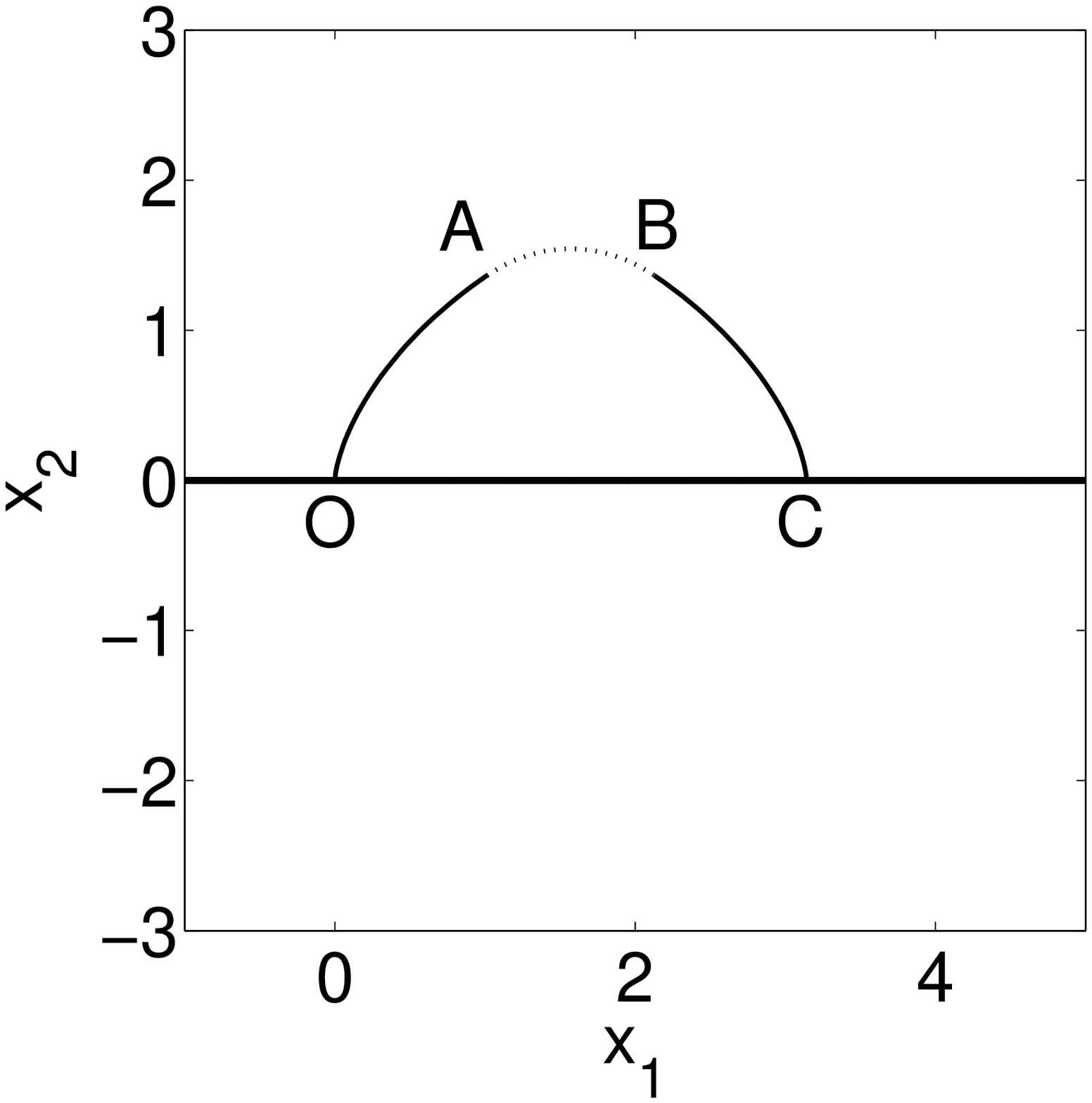}} \\
		\end{tabular}
\caption{Projection of the optimal trajectory on the $x_1x_2$-plane for $\gamma=\pi$. The circles with center $(x_3,\pm 1)$ that generate this projection by rolling on the line $x_2=0$ with velocity  $\dot{x}_3=u=\pm1$ are also shown (a) The first part $OA$ is a cycloid (b) The second part $AB$ is a prolate trochoid, since $A$ lies outside the rolling disc (c) The last part $BC$ is symmetric to the first one (d) Total trajectory.}
 \label{fig:example}
\end{figure*}

There is a simple way to visualize the extremal trajectories in two dimensions. It is based on the observation that the projections of these trajectories on the $x_1x_2$-plane are concatenations of \emph{trochoids}. Recall that a trochoid is the locus of a point at some fixed distance from the center of a circle rolling on a fixed line. Indeed, if we set
\begin{equation*}
y_1=x_1-x_3,\quad y_2=x_2\mp1,
\end{equation*}
for $u=\pm1$,
then we find
\begin{equation*}
y_1^2+y_2^2=\mbox{constant},
\end{equation*}
for each time interval where the control is constant, and
\begin{align}
\label{y_harm_osc1}\dot{y}_1&=y_2,\\
\label{y_harm_osc2}\dot{y}_2&=-y_1.
\end{align}
From the last equations we find that the angular velocity of the rolling circle is $\omega_c=1$. The center of the circle is $(x_3,\pm 1)$, so the horizontal velocity is $v_c=|\dot{x}_3|=|u|=1$ and the radius is $R_c=v_c/\omega_c=1$. The circle rolls without slipping on the line $x_2=0$. In Fig. \ref{fig:example} we plot an extremal trajectory with two switchings, along with the rolling circles that generate it. The first part of the trajectory, $OA$ in Fig. \ref{fig:flash1}, is a cycloid, since $y_1^2+y_2^2=1=R_c^2$. The circle generating this part has center $(x_3,1)$ and rolls to the right since $\dot{x}_3=u=1>0$. When the control switches to $u=-1$, the center of the generating circle becomes $(x_3,-1)$ and it moves to the left, Fig. \ref{fig:flash2}, since now $\dot{x}_3=u=-1<0$. The corresponding trajectory part $AB$ is a prolate trochoid, since the moving point lies outside the rolling disc. After the second switching, the center of the generating circle becomes again $(x_3,1)$ and it rolls to the right, Fig. \ref{fig:flash3}, generating again a cycloid $BC$. The total trajectory $OABC$ is shown in Fig. \ref{fig:total_traj}. The trochoids which compose it are synchronized such that the center of the circle and the point $(x_1,x_2)$ arrive simultaneously at the points $(\gamma,1)$ and $(\gamma,0)$, respectively, so the final point in $\mathbb{R}^3$ is $(\gamma,0,\gamma)$. Note that there is a symmetry between the initial and the final part. As we shall see later, this observation is the key for the optimal solution.

In the next proposition, we use the geometric intuition developed above to calculate the system evolution under an extremal input.

\begin{lem}[Main technical point]
\label{evolution}
Let $0=t_0<t_1<t_2<\dots<t_n=T$. The alternating control input
\begin{equation}
\label{extr_control}
u(t)=(-1)^{j-1},\quad t_{j-1}<t<t_j,\quad j=1,\ldots,n
\end{equation}
drives system (\ref{control_system}) from the origin $x(0)=0$ to the point $x(t_n)$ with coordinates
\begin{align}
\label{X1}x_1(t_n)-x_3(t_n) &= -\sin(t_n)+2\sum_{j=1}^{n-1}(-1)^{j-1}\sin(t_{n}-t_j)\\
\label{X2}x_2(t_n)-(-1)^{n-1} &= -\cos(t_{n})+2\sum_{j=1}^{n-1}(-1)^{j-1}\cos(t_{n}-t_j)\\
\label{X3}x_3(t_n) &= \sum_{j=1}^n(-1)^{j-1}(t_j-t_{j-1})
\end{align}
The control $-u(t)$ drives the system to the symmetric point $-x(t_n)$.
\end{lem}

\begin{IEEEproof}
Note first that since $\dot{x}_3=u$, (\ref{X3}) is obvious for the input (\ref{extr_control}). In order to prove (\ref{X1}) and (\ref{X2}), we use the attached to the rolling circles ``moving" coordinates $y=(y_1,y_2)^T$, where
\begin{align}
\label{y1}y_1&=x_1-x_3,\\
\label{y2}y_2&=x_2-(-1)^{j-1},
\end{align}
for $t_{j-1}<t<t_j, j=1,\ldots,n$. Observe that in each time interval, $y_1$ and $y_2$ satisfy the equations (\ref{y_harm_osc1}) and (\ref{y_harm_osc2}) of the harmonic oscillator, so
\begin{equation}
\label{rotation}
y(t_j^-)=R(t_j-t_{j-1})y(t_{j-1}^+),
\end{equation}
where $R(\tau)$ is the rotation matrix
\begin{equation}
\label{rotation_matrix}
R(\tau)
=
\left[
\begin{array}{cc}
\cos\tau & \sin\tau\\
-\sin\tau & \cos\tau
\end{array}
\right].
\end{equation}
When the control switches, there is a discontinuity in $y_2$
\begin{equation}
\label{discont}
y(t_j^+)=y(t_j^-)+(-1)^{j-1}\left[\begin{array}{c}0\\2\end{array}\right], \quad j=1,\ldots,n-1
\end{equation}
expressing the change in the center of the generating circle, $(x_3,\pm1)\rightarrow (x_3,\mp1)$. Note that
\begin{equation}
\label{init_cond}
y(0^+)=\left[\begin{array}{c}0\\-1\end{array}\right]
\end{equation}
from the definition (\ref{y1}), (\ref{y2}) of $y$ and the initial condition $x(0)=0$.

For $n=1$ we find from (\ref{rotation}), (\ref{rotation_matrix}) and (\ref{init_cond})
\begin{equation*}
y(t_1^-)=R(t_1-t_0)y(t_0^+)=R(t_1)y(0^+)=\left[\begin{array}{c}-\sin(t_1)\\-\cos(t_1)\end{array}\right],
\end{equation*}
thus (\ref{X1}) and (\ref{X2}) hold. Now suppose that they hold for $n$ even (odd), so
\begin{equation}
\label{valid_n}
y(t_n^-)=\left[\begin{array}{c}-\sin(t_n)+2\sum_{j=1}^{n-1}(-1)^{j-1}\sin(t_{n}-t_j)\\-\cos(t_{n})+2\sum_{j=1}^{n-1}(-1)^{j-1}\cos(t_{n}-t_j)\end{array}\right].
\end{equation}
But
\begin{equation}
\label{step}
y(t_n^+)=y(t_n^-)+(-1)^{n-1}\left[\begin{array}{c}0\\2\end{array}\right]=y(t_n^-)\mp\left[\begin{array}{c}0\\2\end{array}\right],
\end{equation}
where the minus (plus) sign in (\ref{step}) corresponds to $n$ even (odd), and
\begin{equation}
\label{induction_step}
y(t_{n+1}^-)=R(t_{n+1}-t_{n})y(t_n^+).
\end{equation}
Using (\ref{valid_n}), (\ref{step}) and (\ref{rotation_matrix}) in (\ref{induction_step}) we find
\begin{align*}
&y_1(t_{n+1}^-)\\
&=-\sin(t_{n+1}-t_n)\cos(t_n)-\cos(t_{n+1}-t_n)\sin(t_n)+\\
&2\{\mp\sin(t_{n+1}-t_n)+\sum_{j=1}^{n-1}(-1)^{j-1}[\sin(t_{n+1}-t_n)\cos(t_n-t_j)+\cos(t_{n+1}-t_n)\sin(t_n-t_j)]\}\\
&=-\sin(t_{n+1})+2[(-1)^{n-1}\sin(t_{n+1}-t_n)+\sum_{j=1}^{n-1}(-1)^{j-1}\sin(t_{n+1}-t_j)]\\
&=-\sin(t_{n+1})+2\sum_{j=1}^{n}(-1)^{j-1}\sin(t_{n+1}-t_j).
\end{align*}
and
\begin{align*}
&y_2(t_{n+1}^-)\\
&=\sin(t_{n+1}-t_n)\sin(t_n)-\cos(t_{n+1}-t_n)\cos(t_n)+\\
&2\{\mp\cos(t_{n+1}-t_n)+\sum_{j=1}^{n-1}(-1)^{j-1}[\cos(t_{n+1}-t_n)\cos(t_n-t_j)-\sin(t_{n+1}-t_n)\sin(t_n-t_j)]\}\\
&=-\cos(t_{n+1})+2[(-1)^{n-1}\cos(t_{n+1}-t_n)+\sum_{j=1}^{n-1}(-1)^{j-1}\cos(t_{n+1}-t_j)]\\
&=-\cos(t_{n+1})+2\sum_{j=1}^{n}(-1)^{j-1}\cos(t_{n+1}-t_j).
\end{align*}
The induction step has been proved. To prove the last statement in the lemma we use the variation of constants formula \cite{Brockett70} for the linear system (\ref{control_system}), which for $x(0)=0$ gives
\begin{equation*}
x(t)=\int_0^{t}u(\sigma)e^{A(t-\sigma)}b\,d\sigma.
\end{equation*}
Obviously, the control $-u(t)$ drives the system to the symmetric point $-x(t_n)$.
\end{IEEEproof}

\begin{theorem}[Optimal solution]
\label{solution}
For the final point $(\gamma,0,\gamma)$, with $2(\rho-1)\pi<\gamma<2\rho\pi, \rho=1,2,\ldots$, problem \ref{problem} has a unique optimal solution with $2\rho$ switchings
\begin{equation}
u(t)=(-1)^{j-1},\quad t_{j-1}<t<t_j,\quad j=1,\ldots,2\rho+1,
\end{equation}
where the constant control time intervals $\tau_j=t_j-t_{j-1}$  are such that the initial and final intervals are equal $\tau_1=\tau_{2\rho+1}=\tau$ and are given by the solution of the following transcendental equation
\begin{equation}
\label{trans_eq}
\frac{2\tau+2(\rho-1)\pi-\gamma}{2\rho-1}=2\tan^{-1}\left(\frac{\sin \tau}{2\rho-\cos \tau}\right),
\end{equation}
while the intermediate intervals are
\begin{equation}
\label{t_2k}
\tau_{2k} = \frac{2\tau+2(\rho-1)\pi-\gamma}{2\rho-1},\quad k=1,\ldots,\rho
\end{equation}
and
\begin{equation}
\label{t_2k+1}
\tau_{2k+1} = 2\pi-\tau_{2k+2}=\frac{2\rho\pi-2\tau+\gamma}{2\rho-1},\quad k=1,\ldots,\rho-1.
\end{equation}
The total minimum transfer time is
\begin{equation}
\label{total_time}
t_{2\rho+1}=\frac{4\rho[\tau+(\rho-1)\pi]-\gamma}{2\rho-1}.
\end{equation}
For $\gamma=2\rho\pi$ the optimal control is $u(t)=1$ and $t_{2\rho+1}=2\rho\pi$.
\end{theorem}

\begin{IEEEproof}
We study first the bang-bang extremals. Consider an extremal control of the form (\ref{extr_control}) with $2\rho$ switchings. From lemma \ref{evolution} we have that the final state satisfies the terminal condition $x(t_{2\rho+1})=(\gamma,0,\gamma)^T$ when
\begin{align}
\label{x11}&-\sin(t_{2\rho+1})+2\sum_{j=1}^{2\rho}(-1)^{j-1}\sin(t_{2\rho+1}-t_j)=0,\\
\label{x2}&-\cos(t_{2\rho+1})+2\sum_{j=1}^{2\rho}(-1)^{j-1}\cos(t_{2\rho+1}-t_j)+1=0,\\
\label{x3}&\sum_{j=1}^{2\rho+1}(-1)^{j-1}(t_j-t_{j-1})=\gamma.
\end{align}
If we multiply (\ref{x11}) by $i=\sqrt{-1}$ and add (\ref{x2}) we obtain
\begin{equation}
\label{complex}
-e^{it_{2\rho+1}}+2\sum_{j=1}^{2\rho}(-1)^{j-1}e^{i(t_{2\rho+1}-t_j)}+1=0.
\end{equation}
We express this relation using the constant control time intervals $\tau_k=t_k-t_{k-1}$. Due to the sinusoidal form with period $2\pi$ of the switching function (\ref{switching_function}), for a bang-bang control it is $0<\tau_k<2\pi, k=1,\ldots,2\rho+1$, as well as $\tau_k+\tau_{k+1}=2\pi$ for $k=2,3,\ldots,2\rho-1$ and $\rho\geq 2$. Also $\tau_{2k}$ are equal for $k=1,2,\ldots,\rho$ and $\rho\geq 2$, while $\tau_{2k+1}$ are equal for $k=1,2,\ldots,\rho-1$ and $\rho\geq 3$. Using these relations, the times appearing in (\ref{complex}) can be expressed as follows
\begin{equation}
\label{tot_time}
t_{2\rho+1}=\sum_{k=1}^{2\rho+1}(t_k-t_{k-1})=\sum_{k=1}^{2\rho+1}\tau_k=\tau_1+2(\rho-1)\pi+\tau_{2\rho}+\tau_{2\rho+1}
\end{equation}
and
\begin{equation}
\label{tranc_times}
t_{2\rho+1}-t_j=\sum_{k>j}^{2\rho+1}(t_k-t_{k-1})=\sum_{k>j}^{2\rho+1}\tau_k=\left\{\begin{array}{lc} (2\rho-j-1)\pi+\tau_{2\rho}+\tau_{2\rho+1}, & j\quad\mbox{odd}\\(2\rho-j)\pi+\tau_{2\rho+1}, &  j\quad\mbox{even}\end{array}\right..
\end{equation}
Using (\ref{tot_time}) and (\ref{tranc_times}) in (\ref{complex}) we obtain
\begin{equation*}
-e^{i(\tau_1+\tau_{2\rho}+\tau_{2\rho+1})}+2\rho e^{i(\tau_{2\rho}+\tau_{2\rho+1})}-2\rho e^{i\tau_{2\rho+1}}+1=0,
\end{equation*}
which leads to
\begin{equation}
\label{exponential}
e^{i\tau_{2\rho}}=\frac{2\rho-e^{-i\tau_{2\rho+1}}}{2\rho-e^{i\tau_1}}.
\end{equation}
By taking the absolute value on both sides in the above equation we obtain
\begin{equation*}
\cos{\tau_1}=\cos{\tau_{2\rho+1}}\Rightarrow \tau_1=\tau_{2\rho+1}\quad\mbox{or}\quad\tau_1=2\pi-\tau_{2\rho+1}.
\end{equation*}
Under the second choice, (\ref{exponential}) takes the form $e^{i\tau_{2\rho}}=1$ which has no solution in $(0, 2\pi)$. So
\begin{equation}
\label{init_final}
\tau_1=\tau_{2\rho+1}\triangleq \tau,\quad 0<\tau<2\pi.
\end{equation}
Using (\ref{init_final}), (\ref{exponential}) becomes
\begin{equation}
\label{exponential1}
e^{i\tau_{2\rho}}=e^{i2\phi},\quad\phi=\tan^{-1}\left(\frac{\sin \tau}{2\rho-\cos \tau}\right),
\end{equation}
where the range of $\tan^{-1}$ is taken to be $(-\pi/2,\pi/2)$. For $0<\tau_{2\rho}<2\pi$, (\ref{exponential1}) implies
\begin{equation}
\label{trans}
\left\{\begin{array}{lc} \tau_{2\rho}=2\phi, &  0<\tau<\pi\,(\phi>0)\\\tau_{2\rho}=2\pi+2\phi, &  \pi\leq \tau<2\pi\,(\phi\leq0)\end{array}\right..
\end{equation}
By expressing (\ref{x3}) in terms of $\tau_k$ and using (\ref{init_final}) and the other relations for these time intervals we obtain
\begin{equation*}
\sum_{k=1}^{2\rho+1}(-1)^{k-1}(t_k-t_{k-1})=\sum_{k=1}^{2\rho+1}(-1)^{k-1}\tau_k=2\tau-\rho\tau_{2\rho}+(\rho-1)(2\pi-\tau_{2\rho})=\gamma,
\end{equation*}
so
\begin{equation}
\label{tau2r}
\tau_{2\rho}=\frac{2\tau+2(\rho-1)\pi-\gamma}{2\rho-1}.
\end{equation}
Using (\ref{exponential1}) and (\ref{tau2r}), (\ref{trans}) becomes
\begin{equation}
\label{transendental}
\left\{\begin{array}{lc} f_{\rho}(\tau)=0, &  0<\tau<\pi\\f_{\rho}(\tau)-2\pi=0, &  \pi\leq \tau<2\pi\end{array}\right.,
\end{equation}
where
\begin{equation}
f_{\rho}(\tau)=\frac{2\tau+2(\rho-1)\pi-\gamma}{2\rho-1}-2\tan^{-1}\left(\frac{\sin \tau}{2\rho-\cos \tau}\right).
\end{equation}
It is
\begin{equation}
f_{\rho}'(\tau)=2\left(\frac{1}{2\rho-1}-\frac{2\rho\cos \tau-1}{4\rho^2-4\rho\cos \tau+1}\right)>f_{\rho}'(0)=0,\quad 0<\tau<2\pi.
\end{equation}
Note that the above derivative attains its minimum value when the second fraction in the parenthesis is maximized. This happens for $\cos \tau=1$, which maximizes the numerator and minimizes the denominator. From the above inequality we conclude that $f_{\rho}(\tau)$ is monotonically increasing in the interval $(0,2\pi)$. Since, additionally,
\begin{equation*}
f_{\rho}(0)=\frac{2(\rho-1)\pi-\gamma}{2\rho-1},\quad f_{\rho}(\pi)=\frac{2\rho\pi-\gamma}{2\rho-1}
\end{equation*}
we observe that for $2(\rho-1)\pi<\gamma<2\rho\pi, \rho=1,2,\ldots$, the equation $f_{\rho}(\tau)=0$ has a unique solution in $(0,\pi)$. On the other hand, 
\begin{equation*}
f_{\rho}(\pi)-2\pi=\frac{2(1-\rho)\pi-\gamma}{2\rho-1},\quad f_{\rho}(2\pi)-2\pi=\frac{2(2-\rho)\pi-\gamma}{2\rho-1}.
\end{equation*}
Note that $f_{\rho}(\pi)-2\pi<0$ for $\rho=1,2,\ldots$ and $\gamma>0$, while $f_{\rho}(2\pi)-2\pi<0$ for $\rho=2,3,\ldots$ and $\gamma>0$. Only for $\rho=1$ and $0<\gamma<2\pi$ it is $f_1(2\pi)-2\pi>0$ and then equation $f_1(\tau)-2\pi=0$ has a solution in $(\pi,2\pi)$. Comparing this with the solution of $f_1(\tau)=0$ for $0<\gamma<2\pi$ and using (\ref{tot_time}), (\ref{init_final}) and (\ref{tau2r}), we find that in both cases the total time is $4\tau-\gamma$ so the latter solution, which lies in $(0,\pi)$, corresponds to a shorter path.

Now consider an extremal control of the form (\ref{extr_control}) with $2\rho-1$ switchings. Working as above we find
\begin{equation}
\label{complex_even}
-e^{it_{2\rho}}+2\sum_{j=1}^{2\rho-1}(-1)^{j-1}e^{i(t_{2\rho}-t_j)}-1=0,
\end{equation}
where now
\begin{equation}
\label{tot_time_even}
t_{2\rho}=\sum_{k=1}^{2\rho}(t_k-t_{k-1})=\sum_{k=1}^{2\rho}\tau_k=\tau_1+2(\rho-1)\pi+\tau_{2\rho}
\end{equation}
and
\begin{equation}
\label{tranc_times_even}
t_{2\rho}-t_j=\sum_{k>j}^{2\rho}(t_k-t_{k-1})=\sum_{k>j}^{2\rho}\tau_k=\left\{\begin{array}{lc} (2\rho-j-1)\pi+\tau_{2\rho}, & j\quad\mbox{odd}\\(2\rho-j-2)\pi+\tau_{2\rho-1}+\tau_{2\rho}, &  j\quad\mbox{even}\end{array}\right..
\end{equation}
Using (\ref{tot_time_even}) and (\ref{tranc_times_even}), (\ref{complex_even}) becomes
\begin{equation*}
-e^{i(\tau_1+\tau_{2\rho})}+2\rho e^{i\tau_{2\rho}}-2(\rho-1) e^{i(\tau_{2\rho-1}+\tau_{2\rho})}-1=0,
\end{equation*}
which leads to
\begin{equation}
\label{exponential_even}
e^{i\tau_1}+2(\rho-1) e^{i\tau_{2\rho-1}}+e^{-i\tau_{2\rho}}=2\rho.
\end{equation}
Observe that this equality holds only if the time intervals $\tau_1, \tau_{2\rho}$ and $\tau_{2\rho-1}$ are integer multiples of $2\pi$, which is not the case since they take values in the interval $(0,2\pi)$.

Next we consider the extremal control $-u(t)$, where $u(t)$ is of the form (\ref{extr_control}) with $2\rho$ switchings and $0<\tau_k<2\pi, k=1,\ldots,2\rho+1$. It is not hard to check that (\ref{exponential})-(\ref{trans}) remain valid, but now
\begin{equation*}
\sum_{k=1}^{2\rho+1}(-1)^{k-1}(t_k-t_{k-1})=\sum_{k=1}^{2\rho+1}(-1)^{k-1}\tau_k=2\tau-\rho\tau_{2\rho}+(\rho-1)(2\pi-\tau_{2\rho})=-\gamma,
\end{equation*}
and hence
\begin{equation}
\label{tau2r_neg}
\tau_{2\rho}=\frac{2\tau+2(\rho-1)\pi+\gamma}{2\rho-1}.
\end{equation}
Using (\ref{exponential1}) and (\ref{tau2r_neg}), (\ref{trans}) becomes
\begin{equation}
\label{transendental_neg}
\left\{\begin{array}{lc} g_{\rho}(\tau)=0, &  0<\tau<\pi\\g_{\rho}(\tau)-2\pi=0, &  \pi\leq \tau<2\pi\end{array}\right.,
\end{equation}
where
\begin{equation}
g_{\rho}(\tau)=\frac{2\tau+2(\rho-1)\pi+\gamma}{2\rho-1}-2\tan^{-1}\left(\frac{\sin \tau}{2\rho-\cos \tau}\right).
\end{equation}
Since
\begin{equation}
g_{\rho}'(\tau)=f_{\rho}'(\tau)>0,\quad 0<\tau<2\pi,
\end{equation}
$g_{\rho}(\tau)$ is monotonically increasing in the interval $(0,2\pi)$. However, since $\gamma>0$, it is also
\begin{equation*}
g_{\rho}(0)=\frac{2(\rho-1)\pi+\gamma}{2\rho-1}>0,\quad g_{\rho}(\pi)=\frac{2\rho\pi+\gamma}{2\rho-1}>0
\end{equation*}
and the equation $g_{\rho}(\tau)=0$ has no solution in $(0,\pi)$. On the other hand, if $\tau$ is a solution of the equation $g_{\rho}(\tau)-2\pi=0$, $\pi\leq \tau<2\pi$, then
\begin{align*}
\tau_{2\rho}+\tau_{2\rho+1}&=\frac{2\tau+2(\rho-1)\pi+\gamma}{2\rho-1}+\tau=\frac{(2\rho+1)\tau+2(\rho-1)\pi+\gamma}{2\rho-1}\\
&\geq\frac{(2\rho+1)\pi+2(\rho-1)\pi+\gamma}{2\rho-1}=\frac{(4\rho-1)\pi+\gamma}{2\rho-1}\\
&>\frac{(4\rho-1)\pi}{2\rho-1}>\frac{(4\rho-2)\pi}{2\rho-1}=2\pi,
\end{align*}
and the resulting control does not correspond to an extremal since it should be $\tau_{2\rho}+\tau_{2\rho+1}\leq 2\pi$.

For an extremal control $-u(t)$, where $u(t)$ of the form (\ref{extr_control}), with $2\rho-1$ switchings and $0<\tau_k<2\pi, k=1,\ldots,2\rho$, it is not hard to check that (\ref{exponential_even}) remains valid, so there is no extremal control sequence of this form.

Finally we study the bang extremals. The constant control $u(t)=1$ drives the system to the points $(2\rho\pi,0,2\rho\pi)$ at $t=2\rho\pi, \rho=1,2,\ldots$ It is the only extremal control that achieves this transfer, thus it is time optimal. Note that $u(t)=-1$ drives the system to the points $(-2\rho\pi,0,-2\rho\pi)$ at $t=2\rho\pi, \rho=1,2,\ldots$ The proof of the theorem is now complete. Relations (\ref{t_2k}) and (\ref{t_2k+1}) are derived using (\ref{tau2r}), while (\ref{total_time}) is easily obtained from (\ref{tot_time}) using (\ref{init_final}) and (\ref{tau2r}).
\end{IEEEproof}

\begin{figure*}[!t]
\centerline{\subfigure[Four switchings]{\includegraphics[width=2.5in]{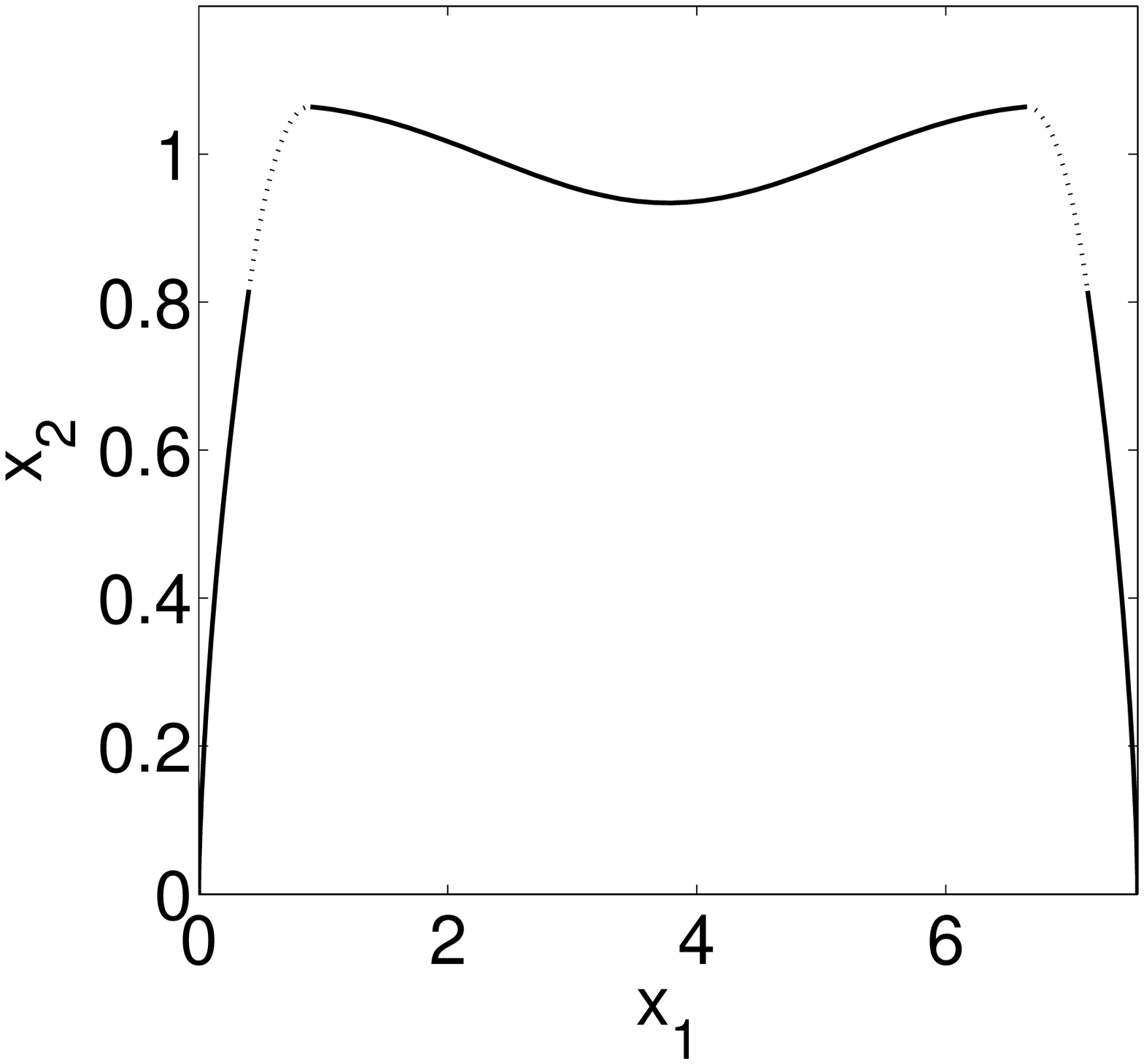}%
\label{fig:Five}}
\hfil
\subfigure[Six switchings]{\includegraphics[width=2.5in]{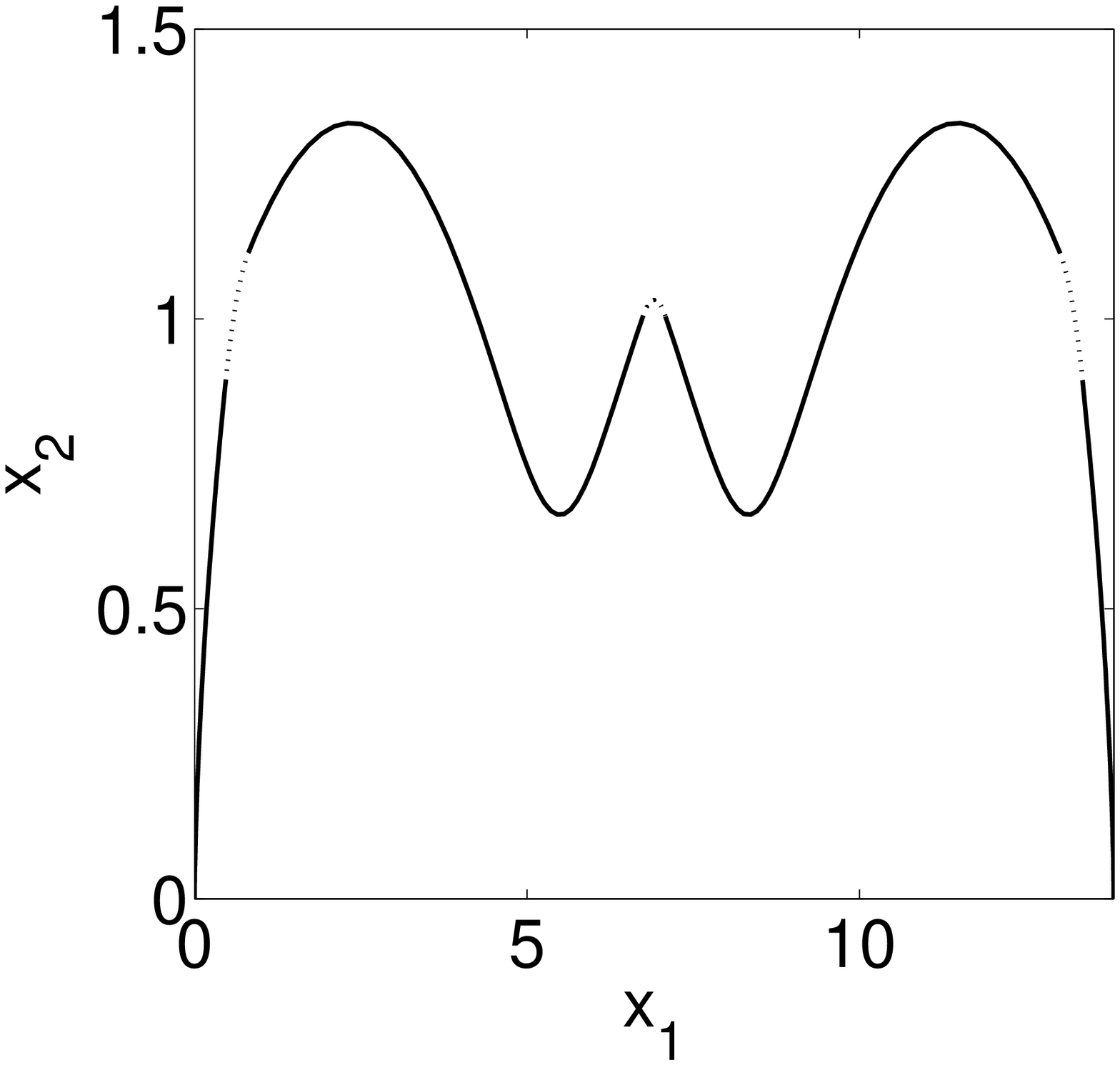}%
\label{fig:Seven}}}
\caption{Projections of the optimal trajectories on the $x_1x_2$-plane for (a) $\gamma=2.4\pi$ and (b) $\gamma=4.4\pi$. Solid (dotted) line corresponds to $u=1$ ($u=-1$).}
\label{fig:five_seven}
\end{figure*}

\begin{figure*}[!t]
\centerline{\subfigure[]{\includegraphics[width=2.5in]{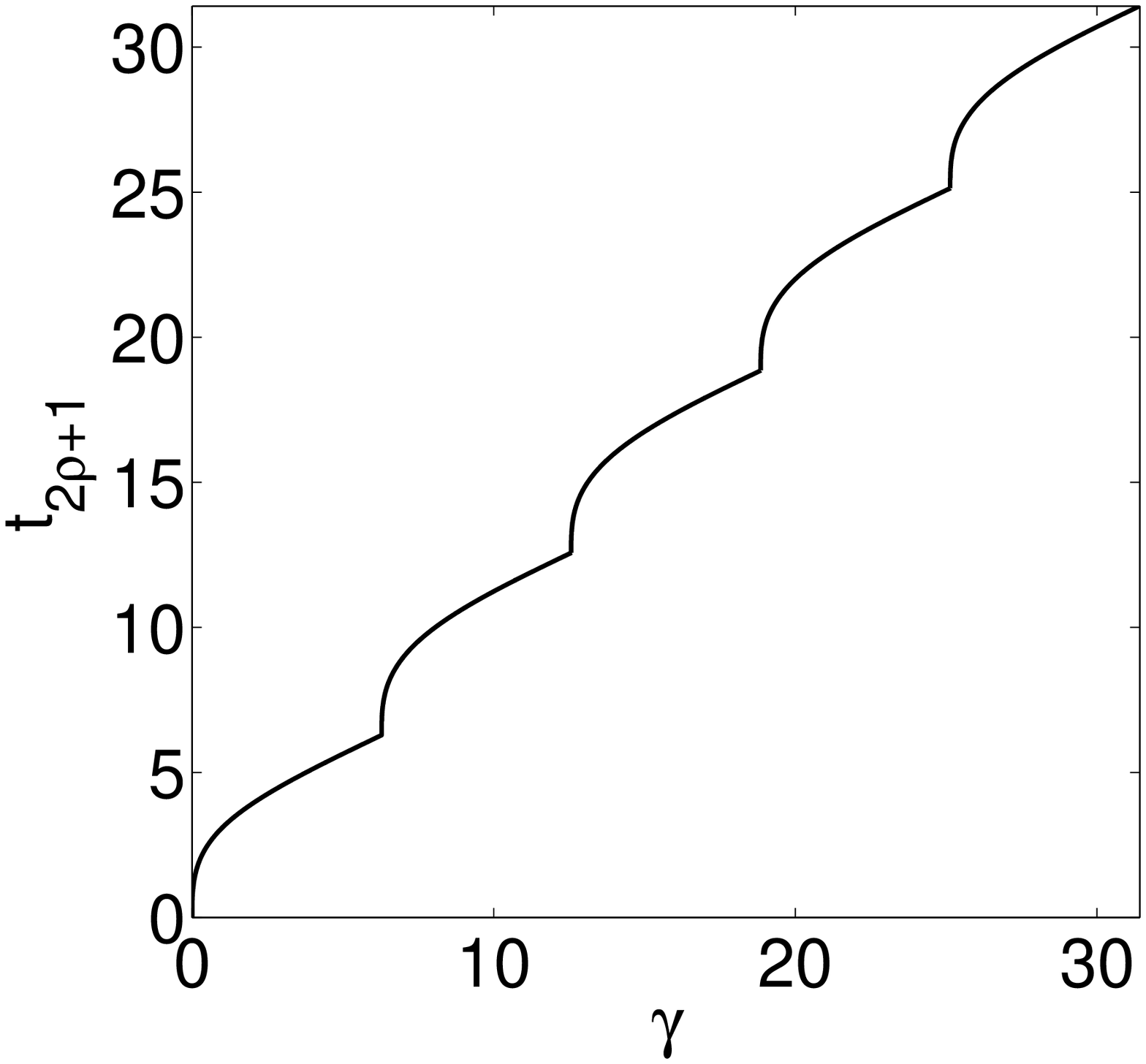}%
\label{fig:time}}
\hfil
\subfigure[]{\includegraphics[width=2.5in]{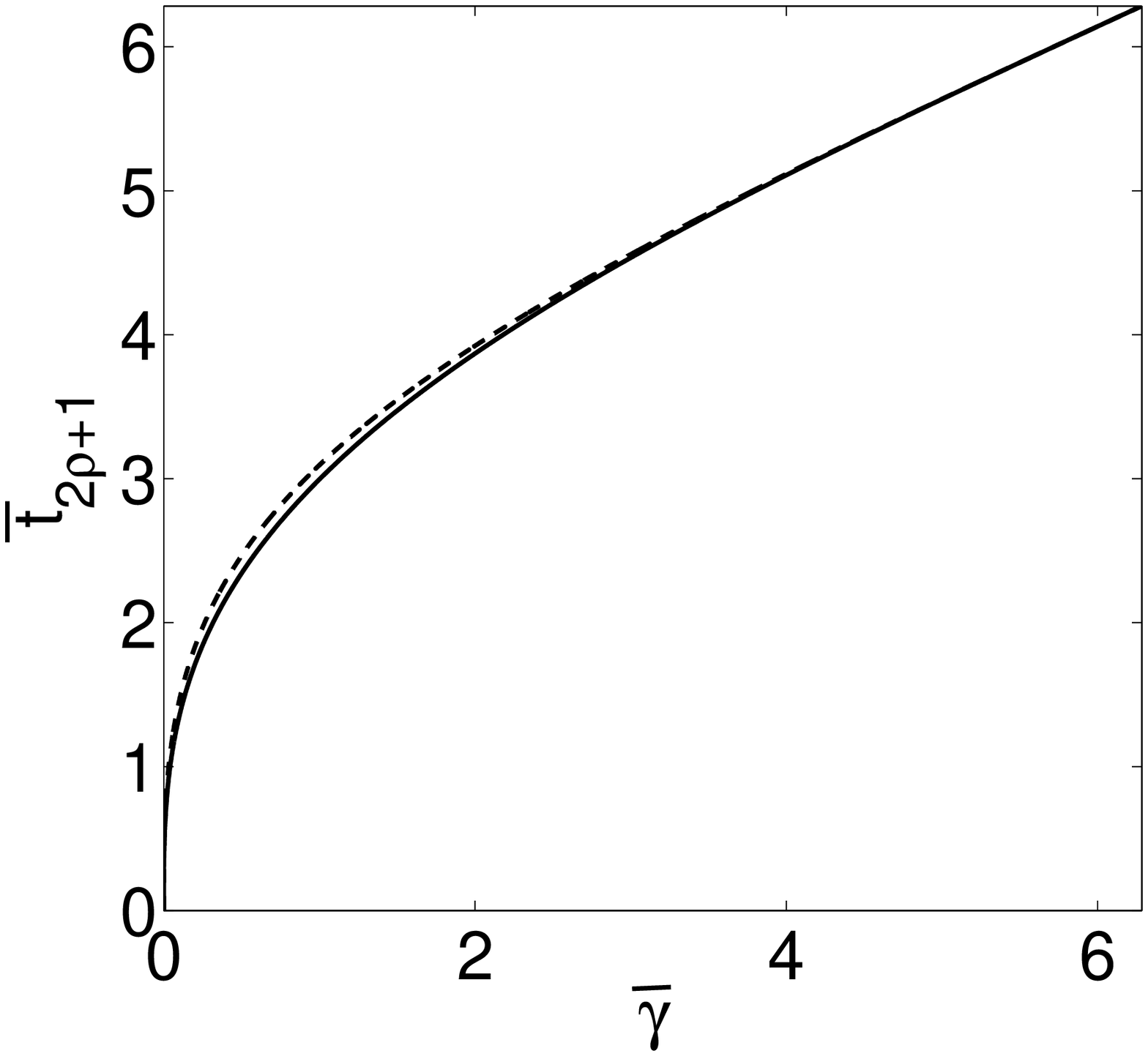}%
\label{fig:time_limit}}}
\caption{(a) Minimum time $t_{2\rho+1}$ as a function of $\gamma\in[0,10\pi]$ (b) $\bar{t}_{2\rho+1}=t_{2\rho+1}-2(\rho-1)\pi$ as a function of $\bar{\gamma}=\gamma-2(\rho-1)\pi$, $\bar{\gamma}\in[0,2\pi]$, where dashed (solid) line corresponds to $\rho=1$ ($\rho\rightarrow\infty$).}
\label{fig:time_and_limit}
\end{figure*}

In Fig. \ref{fig:five_seven} we plot the projections of the optimal trajectories on the $x_1x_2$-plane for $\gamma=2.4\pi$ and $\gamma=4.4\pi$.
In Fig. \ref{fig:time} we plot the minimum time to reach the final point $(\gamma,0,\gamma)$ as a function of $\gamma\in[0,10\pi]$. It is tempting to think that the plot segments for $2(\rho-1)\pi<\gamma<2\rho\pi$ are translations of the initial segment ($0<\gamma<2\pi$) by $2(\rho-1)\pi$ in both axes. But this is not the case, as we shall immediately see. Let
\begin{align}
\bar{\gamma}&=\gamma-2(\rho-1)\pi,\\
\bar{t}_{2\rho+1}&=t_{2\rho+1}-2(\rho-1)\pi.
\end{align}
The transcendental equation (\ref{trans_eq}) can be expressed as
\begin{equation}
\label{trans_bar}
\frac{\tau-\bar{\gamma}/2}{2\rho-1}-\tan^{-1}\left(\frac{\sin \tau}{2\rho-\cos \tau}\right)=0,
\end{equation}
while (\ref{total_time}) gives
\begin{equation}
\bar{t}_{2\rho+1}=\frac{4\rho\tau-\bar{\gamma}}{2\rho-1}.
\end{equation}
In the limit $\rho\rightarrow\infty$, (\ref{trans_bar}) becomes
\begin{equation}
\tau-\sin\tau-\bar{\gamma}/2=0,
\end{equation}
and
\begin{equation}
\label{bar_limit}
\bar{t}_{2\rho+1}\rightarrow 2\tau.
\end{equation}
In Fig. \ref{fig:time_limit} we plot $\bar{t}_{2\rho+1}$ as a function of $\bar{\gamma}\in[0,2\pi]$ for $\rho=1$ (dashed line) and for $\rho\rightarrow\infty$ (solid line). The segments in Fig. \ref{fig:time} approach the limiting case as $\rho$ increases.

\section{Conclusion}

In this article, we formulated the problem of efficient transport of a quantum particle trapped in a harmonic potential moving with a bounded speed, as a minimum-time problem on a linear system with bounded input. We completely solved the corresponding optimal control problem, obtaining an interesting bang-bang solution. Similar approach can be followed for the problem of atom stopping or launching. Additional restrictions on the control, reflecting possible experimental limitations, can be incorporated in the current analysis. The complexity of the resulting optimal control problems, which may increase the difficulty of the analytical study, can be overcome by using a powerful state of the art numerical optimization method based on pseudospectral approximations \cite{Li09,Li11, Ruths11}.

The results obtained here are expected to find application in quantum information processing, where quantum transport between operational sites is an indispensable step. They can also be immediately extended to the efficient transport of Bose-Einstein condensates \cite{Torrontegui11}, where the ability to control the condensates is crucial for their potential use as interferometric sensors \cite{HanselNature, HanselPRA01}.

\section*{Acknowledgment}

The authors would like to thank Prof. J. G. Muga, for bringing the problem of quantum transport in their attention, and Prof. H. Schaettler, for numerous valuable comments and especially for pointing out that there must be an elegant proof, based on geometric intuition, of proposition \ref{evolution}.

\begin{IEEEbiography}{Dionisis Stefanatos}
(M'11) received in 2005 the PhD in Engineering Sciences from Harvard University, where he was a co-recipient of the Eli Jury award and a postdoctoral fellow the next year. Currently is a postdoctoral associate in Washington University St. Louis. His research is focused on the study of control systems that arise from physical problems and especially quantum mechanical applications.
\end{IEEEbiography}

\begin{IEEEbiography}{Jr-Shin Li}
(M'06) received his BS and MS degrees from National Taiwan University, and his PhD degree in Applied Mathematics from Harvard University in 2006. He is currently an Assistant Professor in the Department of Electrical and Systems Engineering with a joint appointment in the Division of Biology \& Biomedical Sciences at Washington University in St. Louis.

His research interests lie in the areas of control theory and optimization. His current work is on the control of complex systems with applications ranging from quantum mechanics and neuroscience to bioinformatics. He is a recipient of the NSF Career Award in 2007 as well as the AFOSR Young Investigator Award in 2009.
\end{IEEEbiography}








\begin{thebibliography}{99}

\bibitem{Mabuchi05}
H.~Mabuchi and N.~Khaneja, ``Principles and applications of control in quantum systems", \emph{Int. J. Robust Nonlinear Control}, vol. 15, pp.~647--667, 2005.

\bibitem{Wiseman93}
H.~M. Wiseman and G.~J. Milburn, ``Quantum theory of optical feedback via homodyne detection", \emph{Phys. Rev. Lett.}, vol. 70, pp.~548--551, 1993.

\bibitem{Doherty00}
A.~C. Doherty, S.~Habib, K.~Jacobs, H.~Mabuchi, and S.~M. Tan, ``Quantum feedback control and classical control theory", \emph{Phys. Rev. A}, vol. 62, 012105, 2000.

\bibitem{Wiseman94}
H.~M. Wiseman and G.~J. Milburn, ``All-optical versus electro-optical quantum-limited feedback", \emph{Phys. Rev. A}, vol. 49, pp.~4110--4125, 1994.

\bibitem{Lloyd00}
S.~Lloyd, ``Coherent quantum feedback", \emph{Phys. Rev. A}, vol. 62, 022108, 2000.

\bibitem{James08}
M.~R. James, H.~I. Nurdin, and I.~R. Petersen, ``$H^\infty$ control of linear quantum stochastic systems", \emph{IEEE Trans. Autom. Control}, vol. 53, no. 8, pp.~1787--1803, Sep. 2008.

\bibitem{Huang83}
G.~M. Huang, T.~J. Tarn, and J.~W. Clark, ``On the controllability of quantum-mechanical systems", \emph{J. Math. Phys.}, vol. 24, pp.~2608--2618, 1983.



\bibitem{Li09IEEE}
J.-S.~Li and N.~Khaneja, ``Ensemble control of Bloch equations", \emph{IEEE Trans. Autom. Control}, vol. 54, no. 3, pp.~528--536, Mar. 2009.


\bibitem{Bloch10}
A.~M. Bloch, R.~W. Brockett, and C.~Rangan, ``Finite controllability of infinite-dimensional quantum systems, \emph{IEEE Trans. Autom. Control}, vol. 55, no. 8, pp.~ 1797--1805, Aug. 2010.


\bibitem{D'Alessandro01}
D.~D'Alessandro and M.~Dahleh, ``Optimal control of two-level quantum systems", \emph{IEEE Trans. Autom. Control}, vol. 46, no. 6, pp.~866--876, Jun. 2001.

\bibitem{Khaneja01}
N.~Khaneja, R.~Brockett, and S.~J. Glaser, ``Time optimal control in spin systems", \emph{Phys. Rev. A}, vol. 63, 032308, 2001.

\bibitem{Boscain02}
U.~Boscain, G.~Charlot, J.-P.~ Gauthier, S.~Guerin, and H.-R.~ Jauslin, ``Optimal control in laser-induced population transfer for two- and three-level quantum systems", \emph{J. Math. Phys.}, vol. 43, pp.~2107--2132, 2002.

\bibitem{Sklarz04}
S.~E. Sklarz, D.~J. Tannor, and N.~Khaneja, ``Optimal control of quantum dissipative dynamics: Analytic solution for cooling the three-level $\Lambda$ system", \emph{Phys. Rev. A}, vol. 69, 053408, 2004.

\bibitem{Stefanatos04}
D.~Stefanatos, N.~Khaneja, and S.~J. Glaser, ``Optimal control of coupled spins in the presence of longitudinal and transverse relaxation", \emph{Phys. Rev. A}, vol. 69, 022319, 2004.

\bibitem{Stefanatos05}
D.~Stefanatos, S.~J. Glaser, and N.~Khaneja, ``Relaxation-optimized transfer of spin order in Ising spin chains", \emph{Phys. Rev. A}, vol. 72, 062320, 2005.

\bibitem{Boscain06}
U.~Boscain and P.~Mason, ``Time minimal trajectories for a spin-1/2 particle in a magnetic field", \emph{J. Math. Phys.}, vol. 47, 062101, 2006.



\bibitem{BonnardIEEE}
B.~Bonnard, M. Chyba, D. Sugny, ``Time-minimal control of dissipative two-level quantum systems: The generic case", \emph{IEEE Trans. Autom. Control}, vol. 54, no. 11, pp. 2598--2610, Nov. 2009.

\bibitem{Fisher09}
R.~Fisher, H.~Yuan, A.~Sp\"{o}rl, and S.~J. Glaser, ``Time-optimal generation of cluster states", \emph{Phys. Rev. A}, vol. 79, 042304, 2009.

\bibitem{Stefanatos09}
D.~Stefanatos, ``Optimal design of minimum-energy pulses for Bloch equations in the case of dominant transverse relaxation", \emph{Phys. Rev. A}, vol. 80, 045401, 2009.

\bibitem{Lapert10}
M.~Lapert, Y.~Zhang, M.~Braun, S.~J. Glaser, and D.~Sugny, ``Singular extremals for the time-optimal control of dissipative spin 1/2 particles", \emph{Phys. Rev. Lett.}, vol. 104, 083001, 2010.

\bibitem{Bonnard10}
B.~Bonnard, O.~Cots, N.~Shcherbakova, and D.~Sugny, ``The energy minimization problem for two-level dissipative quantum systems", \emph{J. Math. Phys.}, vol. 51, 092705, 2010.

\bibitem{Stefanatos10}
D.~Stefanatos and J.-S.~Li, ``Constrained minimum-energy optimal control of the dissipative Bloch equations", \emph{Syst. Control Lett.}, vol. 59, pp.~601--607, 2010.

\bibitem{Stefanatos11}
D.~Stefanatos, H.~Schaettler, and J.-S.~Li, ``Minimum-time frictionless atom cooling in harmonic traps", \emph{SIAM J. Control Optim.}, pp., 2011.



\bibitem{Peirce88}
A.~Peirce, M.~Dahleh, and H.~Rabitz, ``Optimal control of quantum mechanical systems: Existence, numerical approximations, and applications", \emph{Phys. Rev. A}, vol. 37, pp.~4950--4964, 1988.


\bibitem{Khaneja05}
N.~Khaneja, T.~Reiss, C.~Kehlet, T.~Schulte-Herbr\"{u}ggen, and S.~J. Glaser, ``Optimal control of coupled spin dynamics: Design of NMR pulse sequences by gradient ascent algorithms", \emph{J. Magn. Reson.}, vol. 172, pp.~296--305, 2005.

\bibitem{Li09}
J.-S.~Li, J.~Ruths, and D.~Stefanatos, ``A pseudospectral method for optimal control of open quantum systems", \emph{J. Chem. Phys.}, vol. 131, 164110, 2009.

\bibitem{Schulte10}
T.~Schulte-Herbr\"{u}ggen, S.~J. Glaser, G.~Dirr, and U.~Helmke, ``Gradient flows for optimization in quantum information and quantum dynamics: Foundations and applications", \emph{Rev. Math. Phys.}, vol. 22, pp.~597--667, 2010.

\bibitem{Maximov10}
I.~I. Maximov, J.~Salomon, G.~Turinici, and N.~C. Nielsen, ``A smoothing monotonic convergent optimal control algorithm for nuclear magnetic resonance pulse sequence design",
\emph{J. Chem. Phys.}, vol. 132, 084107, 2010.

\bibitem{Li11}
J.-S.~Li, J.~Ruths, T.-Y.~Yu, H.~Arthanari, and G.~Wagner, ``Optimal pulse design in quantum control: A unified computational method", \emph{Proc. Natl. Acad. Sci. U.S.A.}, vol. 108, no. 5, pp.~1879--1884, 2011.


\bibitem{Ruths11}
J.~Ruths and J.-S.~Li, ``A multidimensional pseudospectral method for optimal control of quantum ensembles", \emph{J. Chem. Phys.}, vol. 134, 044128, 2011.

\bibitem{Doria11}
P.~Doria, T.~Calarco, and S.~Montangero, ``Optimal control technique for many-body quantum dynamics", \emph{Phys. Rev. Lett.}, vol. 106, 190501, 2011.



\bibitem{Cirac04}
J.~I. Cirac and P.~Zoller, ``New frontiers in quantum information with atoms and ions", \emph{Physics Today}, vol. 57, pp.~38-44, 2004.

\bibitem{Masuda10}
S.~Masuda and K.~Nakamura, ``Fast-forward of adiabatic dynamics in quantum mechanics", \emph{Proc. R. Soc. A}, vol. 466, pp.~1135--1154, 2010.

\bibitem{Schulz06}
S.~Schulz, U.~Poschinger, K.~Singer, and F.~Schmidt-Kaler, ``Optimization of segmented linear Paul traps and transport of stored particles", \emph{Fortschr. Phys.}, vol. 54, pp.~648--665, 2006

\bibitem{Reichle06}
R.~Reichle, D.~Leibfried, R.~B. Blakestad, J.~Britton, J.~D. Jost, E.~Knill, C.~Langer, R.~Ozeri, S.~Seidelin, and D.~J. Wineland, ``Transport dynamics of single ions in segmented microstructured Paul trap arrays", \emph{Fortschr. Phys.}, vol. 54, pp.~666--685, 2006

\bibitem{Couvert08}
A.~Couvert, T.~Kawalec, G.~Reinaudi and D.~Guéry-Odelin, ``Optimal transport of ultracold atoms in the non-adiabatic regime", \emph{Europhys. Lett.}, vol. 83, 13001, 2008.

\bibitem{TorronteguiPRA11}
E.~Torrontegui, S.~Ib\'{a}\~{n}ez, X.~Chen, A.~Ruschhaupt, D.~Guéry-Odelin, and J.~G. Muga, ``Fast atomic transport without vibrational heating", \emph{Phys. Rev. A}, vol. 83, 013415, 2011.

\bibitem{Murphy09}
M. Murphy, L. Jiang, N. Khaneja, and T. Calarco, ``High-fidelity fast quantum transport with imperfect controls", \emph{Phys. Rev. A}, vol. 79,  020301(R), 2009.

\bibitem{Hansel01}
W.~H\"{a}nsel, J.~Reichel, P.~Hommelhoff, and T.~W. H\"{a}nsch, ``Magnetic Conveyor Belt for Transporting and Merging Trapped Atom Clouds", \emph{Phys. Rev. Lett.}, vol. 86, pp.~608--611, 2001.






\bibitem{Leach77}
P.~G.~L. Leach, ``Invariants and wavefunctions for some time-dependent harmonic oscillator-type Hamiltonians", \emph{J. Math. Phys.}, vol. 18, pp.~1902--1907, 1977.

\bibitem{Merzbacher98}
E.~Merzbacher, \emph{Quantum Mechanics}, 3rd~ed.\hskip 1em plus
  0.5em minus 0.4em\relax New York: John Wiley and Sons, 1998.

\bibitem{Jurdjevic97}
V.~Jurdjevic, \emph{Geometric Control Theory}.\hskip 1em plus
  0.5em minus 0.4em\relax Cambridge: Cambridge University Press, 1997.

\bibitem{Lee67}
E.~B. Lee and L.~Markus, \emph{Foundations of Optimal Control Theory}.\hskip 1em plus
  0.5em minus 0.4em\relax New York: Wiley, 1967.

\bibitem{Pontryagin}
L.~S. Pontryagin, V.~G. Boltyanskii, R.~V. Gamkrelidze, and
E.~F. Mishchenko, \emph{The Mathematical Theory of Optimal Processes}.\hskip 1em plus
  0.5em minus 0.4em\relax New York: Interscience
Publishers, 1962.

\bibitem{Brockett70}
R.~W. Brockett, \emph{Finite Dimensional Linear Systems}.\hskip 1em plus
  0.5em minus 0.4em\relax New York: John Wiley and Sons, 1970.

\bibitem{Torrontegui11}
E.~Torrontegui, X.~Chen, M.~Modugno, S. Schmidt, A.~Ruschhaupt,  and J.~G. Muga, ``Fast transport of Bose-Einstein condensates", e-print arXiv:1103.2532.

\bibitem{HanselNature}
W.~H\"{a}nsel, P.~Hommelhoff, T.~W. H\"{a}nsch, and J.~Reichel, ``Bose-Einstein condensation on a microelectronic chip", \emph{Nature}, vol. 413, pp.~498-501, 2001.

\bibitem{HanselPRA01}
W.~H\"{a}nsel, J.~Reichel, P.~Hommelhoff, and T.~W. H\"{a}nsch, ``Trapped-atom interferometer in a magnetic microtrap",  \emph{Phys. Rev. A}, vol. 64, 063607, 2001.

\end{thebibliography}
\end{document}